\documentclass{article}
\usepackage{graphicx}
\usepackage{epsfig}

\def\igg#1#2#3#4#5{\begin{figure}[!ht]\begin{center}%
\includegraphics[height=#3\textheight]{#1.eps}\hfil
\includegraphics[height=#3\textheight]{#2.eps}
\caption{#5}\label{#4}%
\end{center}\end{figure}}

\def\iggg#1#2#3#4#5#6{\begin{figure}[!ht]\begin{center}%
\includegraphics[height=#4\textheight]{#1.eps}\hfil
\includegraphics[height=#4\textheight]{#2.eps}\hfil
\includegraphics[height=#4\textheight]{#3.eps}
\caption{#6}\label{#5}%
\end{center}\end{figure}}

\begin{document}

\title{Du--Hwang Characteristic Area: Catch-22}
\author{A.\,O.~Ivanov, A.\,A.~Tuzhilin\\
Moscow State University, Faculty of Mechanics and Mathematics}
\date{}

\maketitle

We are forced to write about the D.\,Z.~Du and F.\,K.~Hwang methods of Steiner ratio Gilbert--Pollak Conjecture solution again by the paper of W.~Wu and J.~Zhong~\cite{WuZhong}.  Recall that after the series of publications of many authors, see references in~\cite{ITAlg}, it becomes clear that the proof suggested by D.\,Z.~Du and F.\,K.~Hwang in 1990 contains serious gaps. Those gaps are due to the concept of so-called characteristic area of a minimal Steiner tree, which was introduced informally. All attempts of many authors to cover these gaps by giving formal definitions failed because the resulting object must have several properties contradicting to each other. But the paper~\cite{WuZhong} meets the same problems again and contains the same mistake, so we have decided to explain two main gaps in this small note.

Recall, see details in~\cite{ITCRCbook}, that a Steiner minimal tree $G$ connecting a finite boundary set $M$ of points in the Euclidean plane is a plain tree, whose edges are straight segments meeting at its vertices by angles of at least $120^\circ$. Moreover,  we can always assume that all the vertices of degree $1$ and $2$ belongs to $M$. Such three is said to be {\em full}, if it does not contain vertices of degree $2$, and $M$ coincides with the set of the vertices of degree $1$. Each Steiner minimal tree can be represented as a union of its so-called full components, i.e., the maximal subtrees each of which is a full Steiner minimal tree connecting the corresponding subset of $M$. These components intersect each other by points from $M$, which are the vertices of degree $2$ or $3$ in $G$.

The idea of the characteristic area of a Steiner minimal tree $G$ connecting a finite set $M$ of points in the plane is as follows: we ``walk'' around $G$ and join consecutive points from $M$ by straight segments. In simple cases when $G$ is a full Steiner tree and $n=|M|\le 5$ it is easy to see that we get a polygon  with vertex set $M$, containing $G$.  This polygon is referred as a characteristic area of $G$. The main idea of the Du--Hwang approach is to consider only those spanning trees which belong to the characteristic area (inner spanning trees), construct a deformation of the boundary set together with the characteristic area that leads to some optimal configuration, and to control the behavior of the minimal inner spanning trees to obtain the resulting estimate.  So, to go any further we need to expand the concept of the characteristic area to the cases of  an arbitrary $n$ and non-full Steiner trees.

The case of a full Steiner tree with an arbitrary number of vertices is discussed by several authors, including D.Z.~Du~\cite{Du}. The resulting conclusions are as follows: it is possible to define an immersed planar polygon glued from a family of triangles that satisfy all necessary  properties of the characteristic area. In general, there are several such immersed polygons, but this non-uniqueness difficulty can be overcome, see a discussion in~\cite{ITIzv} and~\cite{Du}.

The main difficulties appear under consideration of a non-full Steiner tree. To struggle forward through the Du--Hwang scheme we need two properties of the characteristic area:
\begin{itemize}
\item the continuity of the length of minimal inner spanning tree on $M$ with respect to deformations of the set $M$, and
\item the monotonicity  of the characteristic areas, i.e. the characteristic area of a total tree must contain a characteristic area of any its full component.
\end{itemize}

The problems with continuity were mentioned in~\cite{ITLup}, \cite{deWet1}, \cite{ITIzv}, and described in details in~\cite{IKMS}. Namely, it turns out that if  a characteristic area for a non-full Steiner tree is defined as a union of the characteristic areas of its full components, then the length of minimal inner spanning tree is not continuous, see example in Figure~\ref{fig:cd1}. Indeed, in the left part of the Figure~\ref{fig:cd1} the unique characteristic area (an embedded planar polygon) of a full Steiner tree is depicted. The minimal inner spanning tree is shown as the thick (blue) polygonal line. Let us deform the boundary set moving the upper point along the incoming edge of the tree. The tree becomes non-full, the characteristic area changes spasmodically to the union of the characteristic areas of the two full components, and the lowest edge of the thick inner tree becomes non-interior for the new characteristic area, Figure~\ref{fig:cd1}, Right. The new minimal inner spanning tree is also shown as the thick polygonal line. So, the length of the minimal inner spanning changes spasmodically also.

\igg{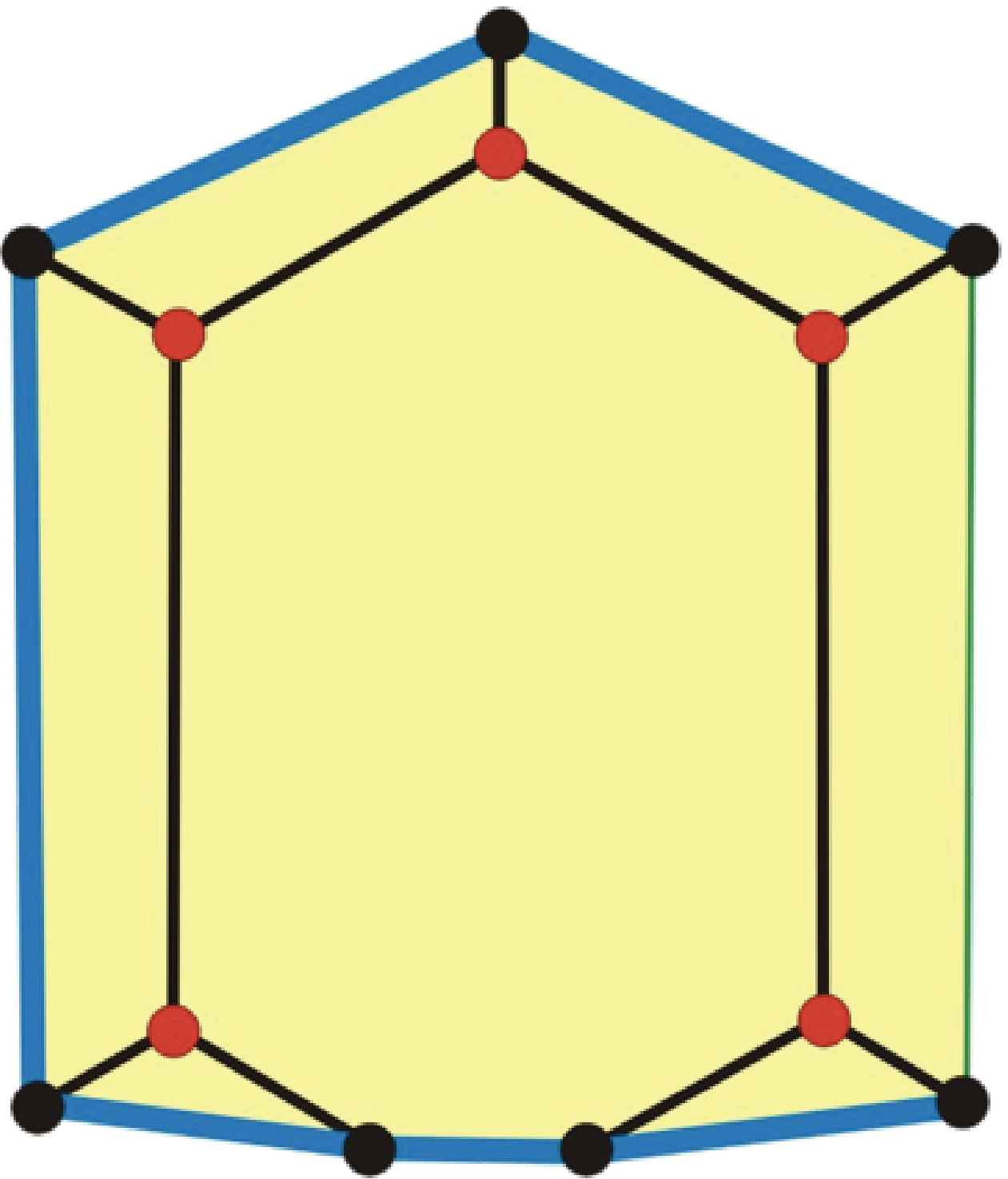}{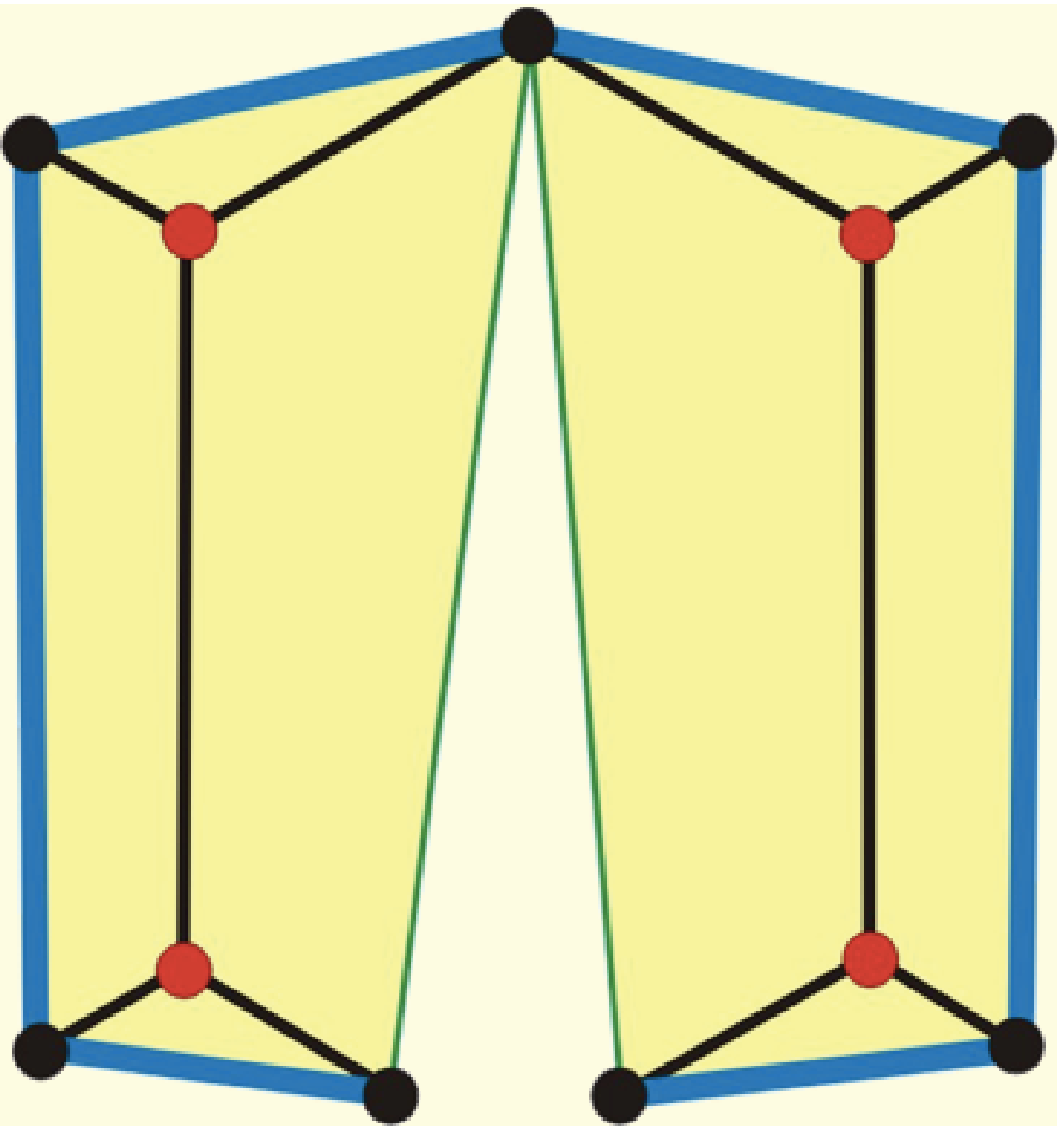}{0.2}{fig:cd1}{{\bf Left:} characteristic area of a full Steiner tree and corresponding minimal inner spanning tree; {\bf Right:} characteristic area for the non-full tree and new minimal inner spanning tree.}

The natural way to overcome this difficulty is to define the characteristic area for non-full Steiner tree as a limit of the characteristic areas of full trees, i.e., by continuity. But this approach leads to non-monotony behavior of characteristic areas. This non-monotonicity effect was mentioned in~\cite{ITIzv}. In Figure~\ref{fig:cd2}, Left, the characteristic area of a full Steiner tree is depicted.  Deform the boundary set as it is shown in Figure~\ref{fig:cd2}, Center (the lowest left vertical edge has been degenerated). The limiting characteristic area is also shown. The resulting  non-full Steiner tree is decomposed into two full components, and the characteristic area of the larger one contains as the characteristic area of the smaller one, so as the characteristic area of the initial non-full Steiner tree, see Figure~\ref{fig:cd2}, Right.

\iggg{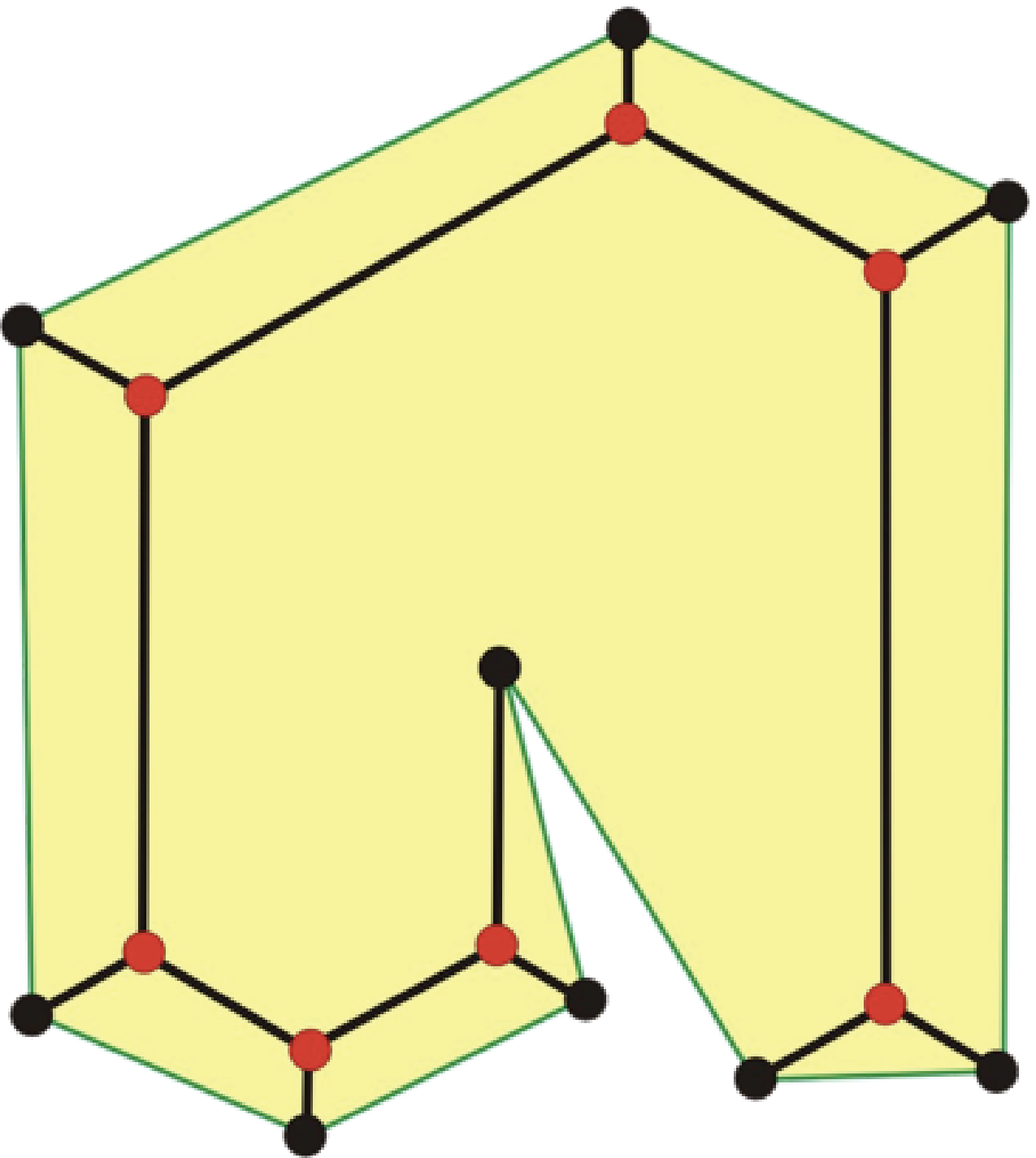}{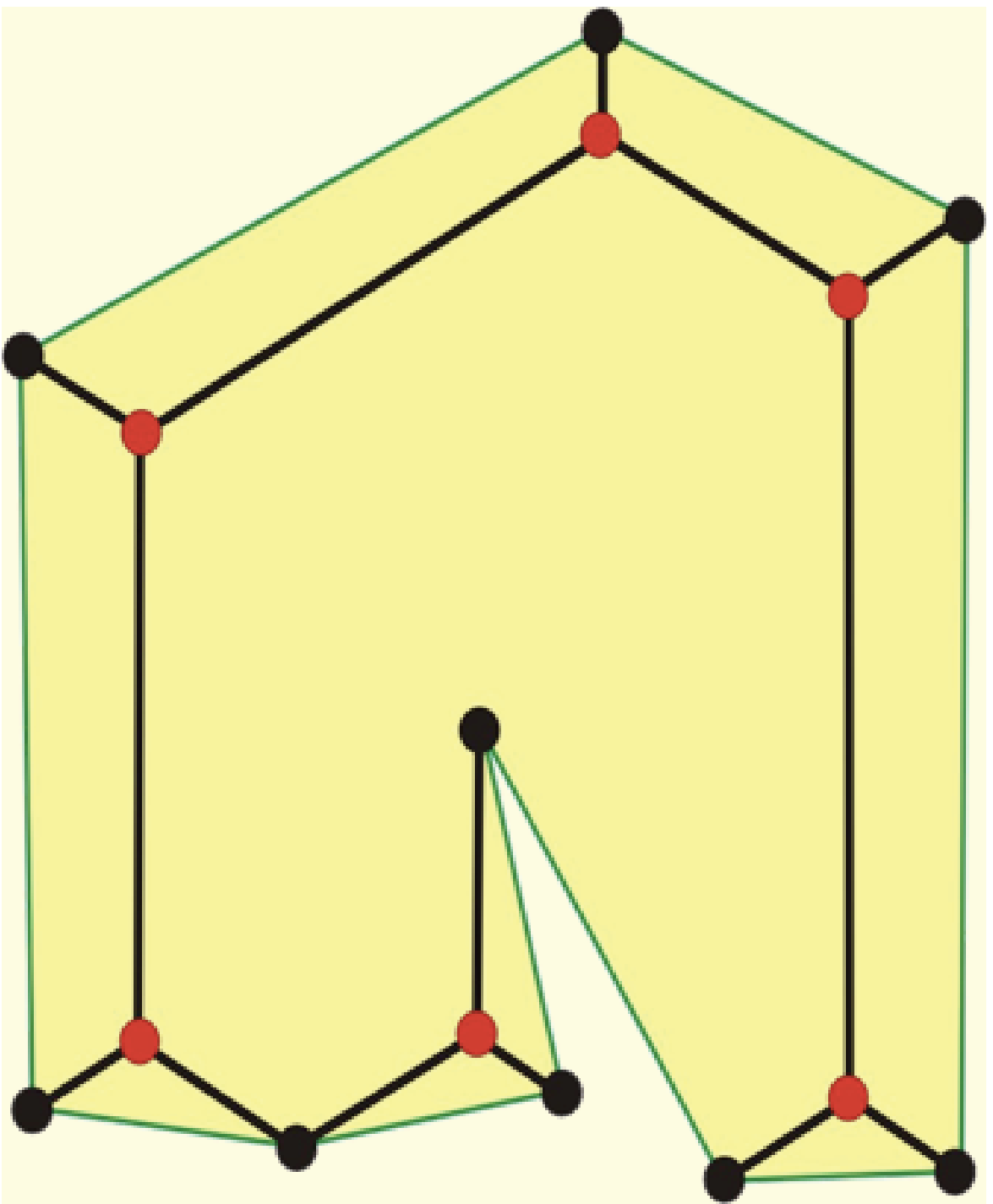}{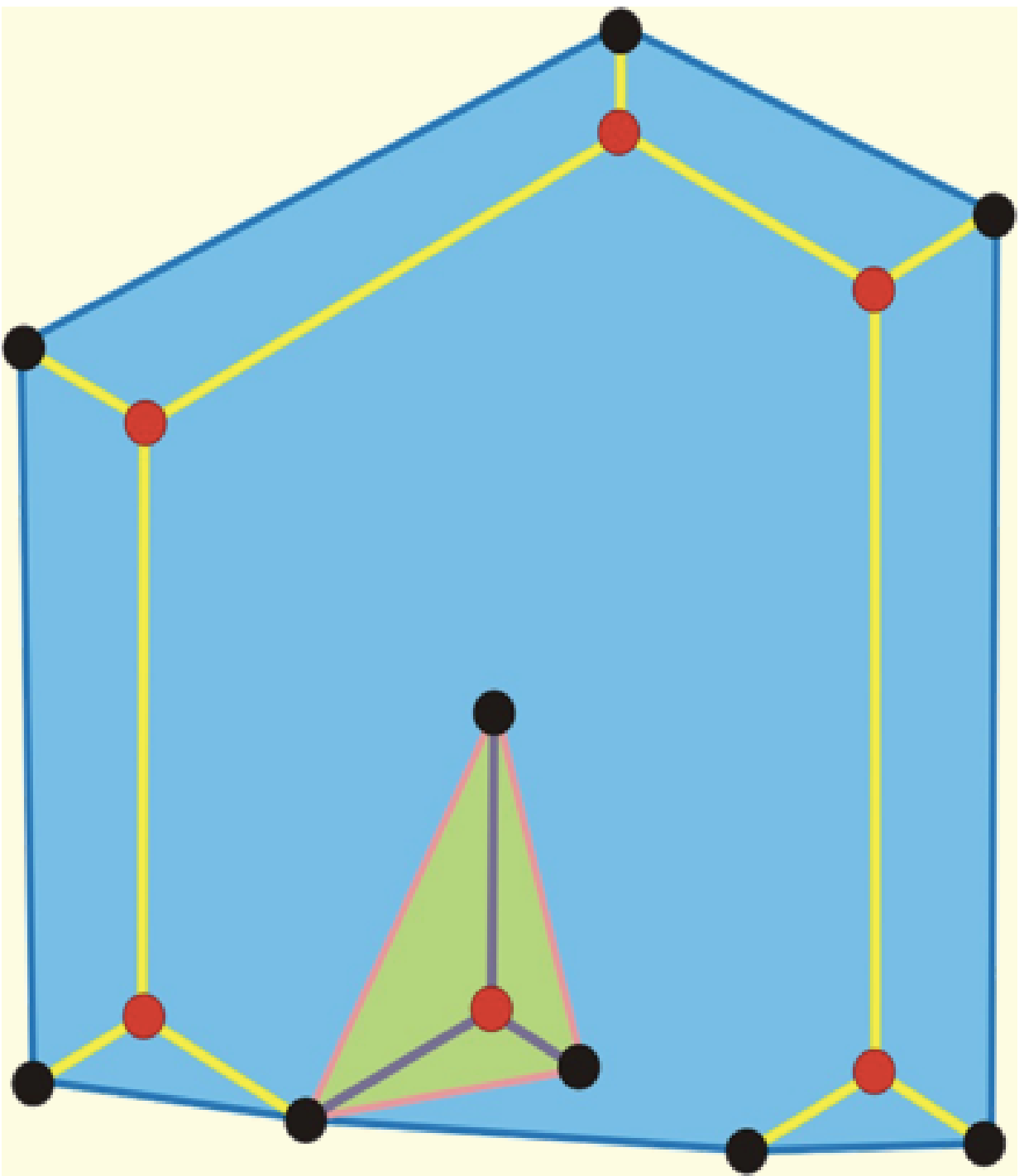}{0.2}{fig:cd2}{{\bf Left:} characteristic area of a full Steiner tree; {\bf Center:} limiting characteristic area; {\bf Right:} characteristic areas of the full components are thick (blue) polygon and small (green) triangle.}

Notice that the monotonicity property is a key one in the Du--Hwang approach since it is necessary to make a reduction of the Steiner ratio estimation to the case of full Steiner minimal trees. W.~Wu and J.~Zhong~\cite{WuZhong} also have used the monotonicity arguments (they need that the union of inner minimal spanning trees be an inner spanning tree for the united boundary sets, Section~3, Case~3), but they just consider an example of a $5$-gon, where the Du--Hwang construction turns out to be valid.

To conclude, let us mention the works~\cite{deWet1}  and~\cite{K} finding the largest number of points in the plane such that the Du--Hwang method works properly. As a result, new proofs of Gilbert--Pollack Conjecture for seven and eight points are obtained. But in general situation we need to find some new ideas and constructions. Some of them are discussed in~\cite{ITGromov}

\end{document}